\documentclass[11pt,a4paper]{article}
\usepackage[T1]{fontenc}
\usepackage{lmodern}
\usepackage{amsmath,amssymb,amsthm,mathtools}
\usepackage[margin=25mm]{geometry}
\usepackage[hidelinks]{hyperref}
\hypersetup{pdftitle={Uniqueness and non-uniqueness of smooth solutions to the incompressible Navier-Stokes equations in critical Besov spaces},pdfsubject={Critical Besov uniqueness and non-uniqueness for Navier-Stokes}}
\usepackage{microtype,needspace}
\usepackage{fancyhdr}
\allowdisplaybreaks[1]
\newtheorem{theorem}{Theorem}

\newtheorem{proposition}{Proposition}[section]
\newtheorem{lemma}[proposition]{Lemma}
\numberwithin{equation}{section}
\theoremstyle{remark}\newtheorem*{remark}{Remark}
\newcommand{\T}{\mathbb T}
\newcommand{\Pj}{\mathbb P}
\newcommand{\Id}{\mathrm{Id}}
\newcommand{\diver}{\operatorname{div}}
\newcommand{\supp}{\operatorname{supp}}
\newcommand{\heat}[1]{e^{#1\Delta}}
\newcommand{\norm}[2]{\left\|#1\right\|_{#2}}
\begin{document}

\title{\begin{center}
{\Large\bfseries Uniqueness and non-uniqueness of smooth solutions\\[3pt]
to the incompressible Navier--Stokes equations\\[3pt]
in critical Besov spaces $\dot B^{-1}_{\infty,q}$}
\end{center}
\hspace{-4mm}
}

\author{Zipeng $\mbox{Chen}^1$ \footnote{Email: chenzp26@mail2.sysu.edu.cn},\quad Song $\mbox{Liu}^1$ \footnote{Email: lius37@mail2.sysu.edu.cn}\quad
	 and\quad
	Zhaoyang $\mbox{Yin}^{1}$\footnote{E-mail: mcsyzy@mail.sysu.edu.cn}\\
    $^1\mbox{School}$ of Science,\\ Shenzhen Campus of Sun Yat-sen University, Shenzhen 518107, China}
 
\date{}
\maketitle
\hrule

\begin{abstract}
We study the incompressible Navier--Stokes equations on
$\mathbb T^d$ with $d\ge2$.
For every $2<q\le\infty$ and every smooth mean-zero background solution $U$,
we construct two distinct mild solutions with the same initial datum
arbitrarily close to $U(0)$ in $\dot B^{-1}_{\infty,q}(\mathbb T^d)$,
both belonging to
$C([0,T];\dot B^{-1}_{\infty,q}(\mathbb T^d))
\cap C^\infty((0,T]\times\mathbb T^d)$.
The time $T>0$ depends on the background solution but is independent
of the prescribed tolerance for the initial data.
This gives non-uniqueness for arbitrarily small initial data
and density of non-unique initial data in the mean-zero divergence-free
subspace of $\dot B^{-1}_{\infty,q}(\mathbb T^d)$ for $2<q<\infty$.
We also prove uniqueness in the same solution class for
$1\le q\le2$, establishing the sharp threshold $q=2$
in this Besov scale.
The uniqueness result extends to
$\dot B^{-1+d/p}_{p,q}(\mathbb T^d)$ for
$1\le p<\infty$ and $1\le q\le\infty$,
under the same continuity and smoothness assumptions.
\end{abstract}

\setcounter{tocdepth}{2}
\tableofcontents
\clearpage

\section{Introduction}

In this paper, we study the incompressible Navier--Stokes equations on $\T^d$:
\begin{equation}\label{eq:NS}
\begin{cases}
\partial_tu+\diver(u\otimes u)-\Delta u+\nabla p=0,\\
\diver u=0,
\end{cases}
\end{equation}
where $d\ge2$, $\T^d=\mathbb R^d/\mathbb Z^d$, and $u=(u_1,\ldots,u_d)$ and $p$ denote the velocity and pressure, respectively.
On $\mathbb R^d$, the equations are invariant under the scaling
\[
u(t,x)\mapsto\lambda u(\lambda^2t,\lambda x),
\qquad p(t,x)\mapsto\lambda^2p(\lambda^2t,\lambda x).
\]
The Besov spaces $\dot B^{-1}_{\infty,q}(\mathbb R^d)$,
$1\le q\le\infty$, are critical for this scaling.

Leray~\cite{Leray} proved the existence of global weak solutions to the three-dimensional Navier--Stokes equations for divergence-free initial data in $L^2(\mathbb R^3)$. These solutions belong to $L^\infty(\mathbb R_+;L^2)\cap L^2(\mathbb R_+;\dot H^1)$ and satisfy the energy inequality
\[
\norm{u(t)}{L^2}^2+2\int_0^t\norm{\nabla u(s)}{L^2}^2\,ds
\le\norm{u(0)}{L^2}^2.
\]
Hopf~\cite{Hopf} established an analogous result on bounded domains with homogeneous Dirichlet boundary conditions. This class of solutions is now known as Leray--Hopf weak solutions.

Convex integration methods for the Euler equations were developed in~\cite{DLS,Isett}. Buckmaster and Vicol~\cite{BV} adapted these methods to construct non-unique finite-energy weak solutions of the three-dimensional Navier--Stokes equations. Their construction uses spatial intermittency to control the viscous term.
For Leray--Hopf weak solutions, the Ladyzhenskaya--Prodi--Serrin criterion~\cite{Ladyzhenskaya,Prodi,Serrin} gives uniqueness if one of the solutions satisfies
\[
u\in L^r(0,T;L^p),
\qquad \frac2r+\frac dp\le1,
\qquad d<p\le\infty.
\]
At $p=d$, the condition $u\in C([0,T];L^d)$ also gives uniqueness; see~\cite{KS}.
Using temporal intermittency, Cheskidov and Luo~\cite{CL22} proved non-uniqueness in $L^r_tL^\infty_x$ for every $1\le r<2$. They also obtained non-uniqueness in $L^\infty_tL^p_x$ for every $1\le p<2$ for the two-dimensional equations~\cite{CL23}.

The mild formulation provides another approach to existence and uniqueness in critical spaces. The relevant spaces satisfy
\[
\dot H^{d/2-1}\hookrightarrow L^d
\hookrightarrow\dot B^{-1+d/p}_{p,\infty}
\hookrightarrow\mathrm{BMO}^{-1}
\hookrightarrow\dot B^{-1}_{\infty,\infty},
\qquad d<p<\infty.
\]
Fujita and Kato~\cite{FK} and Kato~\cite{Kato} established small-data well-posedness in $\dot H^{d/2-1}$ and $L^d$, respectively. Cannone~\cite{Cannone} and Planchon~\cite{Planchon} obtained corresponding results in critical Besov spaces with finite spatial integrability. Koch and Tataru~\cite{KT} proved global well-posedness for sufficiently small initial data in $\mathrm{BMO}^{-1}$.

Ill-posedness also occurs in the critical Besov spaces. Bourgain and Pavlovi\'c~\cite{BP} proved norm inflation in $\dot B^{-1}_{\infty,\infty}$ for the three-dimensional equations. Germain~\cite{Germain} showed that the data-to-solution map is not $C^2$ at the origin in $\dot B^{-1}_{\infty,q}$ for $q>2$, and Yoneda~\cite{Yoneda} proved discontinuity in this range. Wang~\cite{Wang} established norm inflation on $\mathbb R^d$, $d\ge2$, for $1\le q\le2$, although
\[
\dot B^{-1}_{\infty,q}\hookrightarrow\mathrm{BMO}^{-1},
\qquad 1\le q\le2.
\]
Palasek~\cite{Palasek} constructed smooth initial data bounded in $\mathrm{BMO}^{-1}$ whose solutions exhibit arbitrarily large growth in the $\mathrm{BMO}^{-1}$ norm.

Coiculescu and Palasek~\cite{CP} constructed two distinct global solutions of the three-dimensional equations with the same initial datum in $\mathrm{BMO}^{-1}$, both smooth at every positive time. Miao, Nie and Ye~\cite{MNY} obtained a two-dimensional periodic counterpart. Cheskidov, Dai and Palasek~\cite{CDP} constructed solutions that bifurcate from a classical solution through an instantaneous Type~I singularity at a positive time. Fujii~\cite{Fujii} proved non-uniqueness of mild solutions in critical Besov spaces on $\mathbb R^d$, $d\ge3$, without requiring smoothness at every positive time.

A different approach was proposed by Jia and \v Sver\'ak~\cite{JS14,JS15}, based on unstable forward self-similar solutions and a spectral condition for the linearized operator. Albritton, Bru\'e and Colombo~\cite{ABC} proved non-uniqueness of Leray solutions for the forced three-dimensional Navier--Stokes equations. In a recent preprint, Hou, Qin, Song and Zhang~\cite{HQSZ} constructed a family of distinct suitable Leray--Hopf solutions of the unforced three-dimensional equations on $\mathbb R^3$, with a common initial datum in $L^2\cap L^{3,\infty}$ and smoothness at every positive time.

In this paper, we study uniqueness and non-uniqueness in
\[
C([0,T];\dot B^{-1}_{\infty,q}(\T^d))
\cap C^\infty((0,T]\times\T^d).
\]
For $q>2$, we construct two distinct solutions with the same initial datum arbitrarily close to the initial value of a prescribed smooth mean-zero solution. For $1\le q\le2$, we prove uniqueness in this class. The latter result also holds in $\dot B^{-1+d/p}_{p,q}$ for $1\le p<\infty$ and $1\le q\le\infty$.

\subsection{Main results}

\begin{theorem}[Non-uniqueness near a smooth solution]\label{thm:background}
Let $d\ge2$, $2<q_*\le\infty$, and let $U$ be a smooth
mean-zero local solution of \eqref{eq:NS}.
There exists $T=T(U)>0$ such that, for every $\varepsilon>0$,
there exist a mean-zero, divergence-free initial datum $u_0$
and two distinct mild solutions $u^{(1)},u^{(2)}$ on
$[0,T]\times\T^d$ satisfying
\[
u^{(1)}(0)=u^{(2)}(0)=u_0,
\qquad
\norm{u_0-U(0)}{\dot B^{-1}_{\infty,q_*}(\T^d)}
<\varepsilon.
\]
For every $2<q\le\infty$,
\begin{equation}\label{eq:class}
u^{(1)},u^{(2)}\in C([0,T];\dot B^{-1}_{\infty,q}(\T^d))
\cap C^\infty((0,T]\times\T^d).
\end{equation}
Both solutions satisfy
\begin{equation}\label{eq:mild}
u^{(i)}(t)=\heat t u_0-
\int_0^t\heat{(t-s)}\Pj\diver(u^{(i)}\otimes u^{(i)})(s)\,ds,
\qquad i\in\{1,2\},
\end{equation}
where the integral converges in $\dot B^{-1}_{\infty,q}$ for every $2<q\le\infty$.
\end{theorem}

\begin{theorem}[Uniqueness]\label{thm:uniqueness}
Let $d\ge2$, $T>0$, and suppose that
\[
1\le p<\infty,\quad 1\le q\le\infty,
\qquad\text{or}\qquad
p=\infty,\quad 1\le q\le2.
\]
Any two divergence-free mild solutions of \eqref{eq:NS} with the same initial datum and belonging to
\[
C([0,T];\dot B^{-1+d/p}_{p,q}(\T^d))
\cap C^\infty((0,T]\times\T^d)
\]
coincide on $[0,T]$.
\end{theorem}

\begin{remark}
Taking $U=0$ on $[0,1]$ gives non-uniqueness on this interval for arbitrarily small initial data, including smallness in $\dot B^{-1}_{\infty,\infty}$.
For $2<q_*<\infty$, the initial data admitting two distinct solutions in \eqref{eq:class} are dense in the mean-zero divergence-free subspace of $\dot B^{-1}_{\infty,q_*}$. At $q_*=\infty$, this density holds in the $\dot B^{-1}_{\infty,\infty}$ closure of smooth mean-zero divergence-free fields.

The solutions constructed in Theorem~\ref{thm:background} also satisfy
\begin{equation}\label{eq:lorentz-class}
u^{(1)},u^{(2)}\in L^{2,q}(0,T;L^\infty(\T^d)),
\qquad 2<q\le\infty,
\end{equation}
where $L^{2,q}$ is the Lorentz space in time. In particular, they belong to $L^p(0,T;L^\infty)$ for every $0<p<2$.
\end{remark}

\begin{remark}
Theorem~\ref{thm:uniqueness} concerns uniqueness in the stated solution class. It is compatible with Wang's norm-inflation result in $\dot B^{-1}_{\infty,q}$, $1\le q\le2$~\cite{Wang}: the contraction argument takes place in the Koch--Tataru space and does not imply continuous dependence in the Besov norm.
\end{remark}

\subsection{Main ideas}

Following~\cite{CLY}, we combine heat-dominated flows with inverse cascade-dominated flows. For the initial data constructed in Section~\ref{sec:principal-construction}, set
\[
v_{k,j}(t)=e^{-(2\pi N_j)^2t}v_{k,j}^0,
\qquad
\bar v_{k-1}(t)=\frac1{n_k-n_{k-1}}
\sum_{j=n_{k-1}+1}^{n_k}e^{-2(2\pi N_j)^2t}v_{k-1}^0,
\]
for $n_{k-1}<j\le n_k$, with $v_k=\sum_{j=n_{k-1}+1}^{n_k}v_{k,j}$.
The geometric decomposition ensures that the leading nonoscillatory part of $\Pj\diver(v_k\otimes v_k)$ is
\[
\frac2{n_k-n_{k-1}}\sum_{j=n_{k-1}+1}^{n_k}
(2\pi N_j)^2e^{-2(2\pi N_j)^2t}v_{k-1}^0
=-\partial_t\bar v_{k-1}.
\]
We therefore define the principal parts by
\[
v^{(1)}=\sum_{\substack{k\ge0\\k\ \mathrm{even}}}v_k
+\sum_{\substack{k\ge1\\k\ \mathrm{odd}}}\bar v_k,
\qquad
v^{(2)}=\sum_{\substack{k\ge1\\k\ \mathrm{odd}}}v_k
+\sum_{\substack{k\ge0\\k\ \mathrm{even}}}\bar v_k.
\]
At $t=0$,
\[
v_k(0)=\bar v_k(0)=v_k^0,
\qquad v^{(1)}(0)=v^{(2)}(0)=\sum_{k\ge0}v_k^0=V^0.
\]

The main difference from~\cite{CLY} is the construction of the initial data. At each step, we distribute the interaction over $n_k-n_{k-1}$ widely separated frequencies. This gives
\[
\norm{v_k^0}{\dot B^{-1}_{\infty,q}}
\lesssim N_{n_{k-1}}^{1/2}(n_k-n_{k-1})^{1/q-1/2},
\qquad k\ge1.
\]
For $q>2$, choosing the layer endpoints sufficiently large gives
\[
V^0\in\bigcap_{2<q\le\infty}\dot B^{-1}_{\infty,q},
\qquad \norm{V^0}{\dot B^{-1}_{\infty,q_*}}\lesssim A^{-1}.
\]
This replaces the additional frequency decay used in the fractional construction.

We estimate the residual defined by
\[
(\partial_t-\Delta)v^{(i)}+\Pj\diver(v^{(i)}\otimes v^{(i)})=-\mathcal F^{(i)}.
\]
Propositions~\ref{prop:construction} and~\ref{prop:residual} give
\[
\begin{gathered}
\left\|\int_0^te^{(t-s)\Delta}\mathcal F^{(i)}(s)\,ds\right\|_{X_1}
+\norm{\mathcal F^{(i)}}{L^1(0,1;\dot B^{-1}_{\infty,2})}\le A^{-1},\\
\sup_{0<t\le1}\sqrt t\,\norm{v^{(i)}(t)}{L^\infty}\lesssim A^{-1},
\end{gathered}
\]
where $\norm\omega{X_T}=\sup_{0<t\le T}t^{1/4}\norm{\omega(t)}{L^\infty}$.
Using these bounds, we construct perturbations of a smooth background solution $U$ by a fixed-point argument:
\[
u^{(i)}=U+v^{(i)}+\omega^{(i)},
\qquad \omega^{(i)}(0)=0,
\qquad \norm{\omega^{(i)}}{X_T}\le4A^{-1}.
\]
Here $T>0$ depends only on $U$ and $T_0$. Taking $A$ sufficiently large gives
\[
\norm{u_0-U(0)}{\dot B^{-1}_{\infty,q_*}}
=\norm{V^0}{\dot B^{-1}_{\infty,q_*}}<\varepsilon,
\qquad t_0=(2\pi A)^{-2}<T.
\]
The seed flow separates the solutions at $t_0$:
\[
\norm{u^{(1)}(t_0)-u^{(2)}(t_0)}{L^\infty}
\ge e^{-1}-8\sqrt{2\pi}\,A^{-1/2}>0.
\]
Approximating the initial data by smooth fields also gives the density assertion.

For uniqueness, we first reduce to mean-zero solutions. The strong initial trace and positive-time smoothness then give
\[
\lim_{h\downarrow0}\norm{e^{t\Delta}u_0}{\mathrm{KT}(h)}=0,
\]
where the local Koch--Tataru norm is defined in Section~\ref{sec:uniqueness}.
We apply the fixed-point theorem to the solutions restarted at $\delta>0$. The difference estimate in the same space, together with $u^{(i)}(\delta)\to u_0$, gives uniqueness as $\delta\downarrow0$.

The paper is organized as follows. Section~\ref{sec:principal-construction} constructs the initial data and principal parts and estimates the residuals. Section~\ref{sec:correction} constructs the perturbations. The main results are proved in Section~\ref{sec:proofs}. Appendices~\ref{app:geometry} and~\ref{app:analysis} contain the geometric identities and the Fourier and heat estimates used in the proofs.

\paragraph{Notation.}
We use the Fourier convention $e^{2\pi i n\cdot x}$. Let $(\dot\Delta_\ell)_{\ell\in\mathbb Z}$ be a homogeneous Littlewood--Paley decomposition with Fourier supports in $\{\tfrac34 2^\ell\le|n|\le\tfrac83 2^\ell\}$. For a periodic distribution $f$, $s\in\mathbb R$ and $1\le p,q\le\infty$, set
\[
\norm f{\dot B^s_{p,q}}
=\left\|\bigl(2^{s\ell}\norm{\dot\Delta_\ell f}{L^p}\bigr)_{\ell\in\mathbb Z}\right\|_{\ell^q}.
\]
On $\T^d$,
\[
\dot\Delta_\ell=0\quad(\ell\le-2),
\qquad f-\int_{\T^d}f=\sum_{\ell\ge-1}\dot\Delta_\ell f.
\]
For mean-zero $f$, the homogeneous and inhomogeneous norms are equivalent; in particular,
\[
\norm f{\dot B^{-1}_{\infty,q}}\le C_q\norm f{L^\infty}.
\]
The Leray projection has symbol
\[
\Pj(n)=\Id-\frac{n\otimes n}{|n|^2}\quad(n\ne0),\qquad \Pj(0)=\Id.
\]
The inverse Laplacian acts on mean-zero distributions. We use $L^{2,q}(0,T;L^\infty)$ for the Lorentz space in time applied to $t\mapsto\norm{u(t)}{L^\infty}$. Constants may depend on the fixed exponents and dimension, but not on the parameters $A$, $k$ or $j$.

\section{Construction of the principal parts}\label{sec:principal-construction}

We carry out the construction on $\T^2$. The extension $(v_1,v_2,0,\ldots,0)$, independent of $x_3,\ldots,x_d$, commutes with $\Delta$ and $\Pj$ and preserves the Besov estimates up to dimensional constants.

\subsection{Initial data}

Denote $S(f)=\nabla f+(\nabla f)^{\mathsf T}$ and $\nabla^\perp=(-\partial_2,\partial_1)$. For mean-zero divergence-free $f$,
\[
\diver S(\Delta^{-1}f)=f.
\]
Let $\Lambda\subset\mathbb S^1\cap\tfrac15\mathbb Z^2$ and $a_\xi$ be as in Lemma~\ref{lem:geometry}, with $\bar\xi=(-\xi_2,\xi_1)$. They satisfy $\Lambda=-\Lambda$, $a_{-\xi}=a_\xi$ and
\begin{equation}\label{eq:geometry}
\sum_{\xi\in\Lambda}a_\xi(R)^2\bar\xi\otimes\bar\xi=R,
\qquad \max_{r,s}|R_{rs}-\Id_{rs}|\le1/10.
\end{equation}

Fix $2<q_*<\infty$ and set $m=14$. Choose $A\in5\mathbb N$ sufficiently large and define $N_j=A^{320^j}$ for $j\ge0$. With $n_0=0$, set
\begin{equation}\label{eq:layer-endpoints}
n_k=n_{k-1}+
\left\lceil\bigl(2^kA N_{n_{k-1}}^{1/2}\bigr)^{2+\frac{4k}{q_*-2}}\right\rceil,
\qquad k\ge1.
\end{equation}
Since $A\in5\mathbb N$ and $\Lambda\subset\tfrac15\mathbb Z^2$, the frequencies $N_j\xi$ belong to $\mathbb Z^2$.

Set
\begin{equation}\label{eq:seed}
v_0^0=\frac1{2\pi N_0}\nabla^\perp
\left(\sum_{\xi\in\Lambda}a_\xi(\Id)e^{2\pi iN_0\xi\cdot x}\right).
\end{equation}
Given $v_{k-1}^0$, for $n_{k-1}<j\le n_k$ set
\begin{equation}\label{eq:carrier}
v_{k,j}^0
=\frac{\sqrt{2}\,N_{n_{k-1}}^{1/2}}{\sqrt{n_k-n_{k-1}}}\,
\nabla^\perp\Bigg[\sum_{\xi\in\Lambda}
 a_\xi\!\left(\Id+\frac{S(\Delta^{-1}v_{k-1}^0)}{N_{n_{k-1}}}\right)
 e^{2\pi iN_j\xi\cdot x}\Bigg],
\end{equation}
and set $v_k^0=\sum_{j=n_{k-1}+1}^{n_k}v_{k,j}^0$.

\begin{proposition}[Initial data]\label{prop:initial}
For sufficiently large $A$, \eqref{eq:seed}--\eqref{eq:carrier} define smooth, real-valued, mean-zero, divergence-free fields satisfying
\begin{equation}\label{eq:carrier-bound}
\norm{v_k^0}{\dot B^{-1}_{\infty,q}}
\le C_qN_{n_{k-1}}^{1/2}(n_k-n_{k-1})^{1/q-1/2},
\qquad k\ge1,\quad 1\le q\le\infty.
\end{equation}
Moreover,
\begin{equation}\label{eq:V0}
V^0=\sum_{k\ge0}v_k^0
\end{equation}
converges in $\dot B^{-1}_{\infty,q}$ for every $2<q\le\infty$, with
\begin{equation*}
\norm{V^0}{\dot B^{-1}_{\infty,q_*}}\le C_{q_*}A^{-1}.
\end{equation*}
\end{proposition}

\begin{proof}
We prove by induction that the construction is well defined and, for $k\ge1$,
\begin{equation}\label{eq:amplitude-induction}
\begin{aligned}
\norm{S(\Delta^{-1}v_{k-1}^0)}{L^\infty}
&\le \tfrac1{10}N_{n_{k-1}},\\
\left\|a_\xi\!\left(
\Id+\frac{S(\Delta^{-1}v_{k-1}^0)}{N_{n_{k-1}}}
\right)\right\|_{H^s}
&\le N_{n_{k-1}}^s,
\qquad \xi\in\Lambda,\quad s=m,m+1.
\end{aligned}
\end{equation}
For $k=1$, \eqref{eq:seed} and $N_0=A$ give
\[
\left\|\frac{S(\Delta^{-1}v_0^0)}{A}\right\|_{L^\infty}
\le CA^{-2},
\qquad
\left\|\frac{S(\Delta^{-1}v_0^0)}{A}\right\|_{H^s}
\le C_sA^{s-2},
\qquad s=m,m+1.
\]
Applying the composition estimate, we obtain
\[
\left\|a_\xi\!\left(
\Id+\frac{S(\Delta^{-1}v_0^0)}{A}
\right)\right\|_{H^s}
\le C_s(1+A^{s-2})\le A^s,
\qquad \xi\in\Lambda,\quad s=m,m+1.
\]
Assume that \eqref{eq:amplitude-induction} holds at step $k$. Applying the product rule to \eqref{eq:carrier}, we obtain, for $n_{k-1}<j\le n_k$,
\begin{equation}\label{eq:carrier-infinity}
\norm{v_{k,j}^0}{L^\infty}
\le C\frac{N_{n_{k-1}}^{1/2}}{\sqrt{n_k-n_{k-1}}}
\bigl(N_j+N_{n_{k-1}}^m\bigr)
\le C\frac{N_{n_{k-1}}^{1/2}}{\sqrt{n_k-n_{k-1}}}N_j.
\end{equation}
Similarly, for $r=m-1,m$,
\begin{equation}\label{eq:spatial-bounds}
\begin{aligned}
\norm{v_{k,j}^0}{H^r}
&\le C\frac{N_{n_{k-1}}^{1/2}}{\sqrt{n_k-n_{k-1}}}N_j^{r+1},\\
\norm{v_k^0}{H^r}
&\le C\frac{N_{n_{k-1}}^{1/2}}{\sqrt{n_k-n_{k-1}}}N_{n_k}^{r+1}.
\end{aligned}
\end{equation}

Combining \eqref{eq:inverse-bounds}, \eqref{eq:carrier-infinity} and \eqref{eq:layer-endpoints}, we obtain
\[
\frac{\norm{S(\Delta^{-1}v_k^0)}{L^\infty}}{N_{n_k}}
\le C\frac{\norm{v_k^0}{L^\infty}}{N_{n_k}}
\le C\frac{N_{n_{k-1}}^{1/2}}{\sqrt{n_k-n_{k-1}}}
\le C2^{-k}A^{-1}\le\tfrac1{10}.
\]
For $s=m,m+1$, \eqref{eq:inverse-bounds} and \eqref{eq:spatial-bounds} give
\begin{equation*}
\left\|\frac{S(\Delta^{-1}v_k^0)}{N_{n_k}}\right\|_{H^s}
\le C_s\frac{N_{n_{k-1}}^{1/2}}{\sqrt{n_k-n_{k-1}}}
\frac{\sum_{j=n_{k-1}+1}^{n_k}N_j^s}{N_{n_k}}
\le C_sN_{n_k}^{s-1}.
\end{equation*}
Applying the composition estimate again gives
\[
\left\|a_\xi\!\left(\Id+
\frac{S(\Delta^{-1}v_k^0)}{N_{n_k}}\right)\right\|_{H^s}
\le C_s(1+N_{n_k}^{s-1})\le N_{n_k}^s.
\]

This completes the induction. The fields are real-valued by $a_{-\xi}=a_\xi$, and mean-zero and divergence-free by \eqref{eq:seed}--\eqref{eq:carrier}.

We next prove \eqref{eq:carrier-bound}. For the seed,
\[
\norm{v_0^0}{L^\infty}\le C,
\qquad
\norm{v_0^0}{\dot B^{-1}_{\infty,q}}\le C_qA^{-1},
\qquad 1\le q\le\infty.
\]
For $k\ge1$ and $n_{k-1}<j\le n_k$, \eqref{eq:carrier} gives
\[
\widehat{v_{k,j}^0}(n)
=\frac{2\sqrt2\,\pi iN_{n_{k-1}}^{1/2}}{\sqrt{n_k-n_{k-1}}}
(-n_2,n_1)\sum_{\xi\in\Lambda}
\widehat{a_\xi\!\left(\Id+
\frac{S(\Delta^{-1}v_{k-1}^0)}{N_{n_{k-1}}}\right)}(n-N_j\xi).
\]
For $|n|\notin(N_j/2,2N_j)$ and $\xi\in\Lambda$,
we have $|n-N_j\xi|\ge N_j/2$ and $|n|\le3|n-N_j\xi|$.
By \eqref{eq:fourier-tail} and $N_{n_{k-1}}\le N_j^{1/320}$,
\begin{equation*}
\sum_{|n|\notin(N_j/2,2N_j)}
|\widehat{v_{k,j}^0}(n)|
\le C\frac{N_{n_{k-1}}^{1/2}}{\sqrt{n_k-n_{k-1}}}
N_{n_{k-1}}^mN_j^{2-m}
\le C\frac{N_{n_{k-1}}^{1/2}}{\sqrt{n_k-n_{k-1}}}N_j^{-1}.
\end{equation*}
Applying a smooth cutoff equal to one on $N_j/2\le|n|\le2N_j$ and supported in $N_j/4\le|n|\le4N_j$, we obtain
\[
\begin{aligned}
\norm{v_{k,j}^0}{\dot B^{-1}_{\infty,q}}
&\le
C_qN_j^{-1}\norm{v_{k,j}^0}{L^\infty}
+C_q\sum_{|n|\notin(N_j/2,2N_j)}
|\widehat{v_{k,j}^0}(n)|\\
&\le
C_q\frac{N_{n_{k-1}}^{1/2}}{\sqrt{n_k-n_{k-1}}}.
\end{aligned}
\]
Since $N_{j+1}=N_j^{320}$, the Littlewood--Paley decomposition gives,
for scalars $(c_j)_{n_{k-1}<j\le n_k}$,
\begin{equation}\label{eq:carrier-multipliers}
\begin{aligned}
&\Bigl\|\sum_{j=n_{k-1}+1}^{n_k}c_jv_{k,j}^0\Bigr\|_{\dot B^{-1}_{\infty,q}}\\
&\quad\le C_q\bigl\|\bigl(c_jN_j^{-1}\|v_{k,j}^0\|_{L^\infty}\bigr)_{j=n_{k-1}+1}^{n_k}\bigr\|_{\ell^q}
+C_q\sum_{j=n_{k-1}+1}^{n_k}|c_j|
\sum_{|n|\notin(N_j/2,2N_j)}|\widehat{v_{k,j}^0}(n)|\\
&\quad\le C_q\frac{N_{n_{k-1}}^{1/2}}{\sqrt{n_k-n_{k-1}}}
\Bigl(\norm{(c_j)}{\ell^q}
+\sum_{j=n_{k-1}+1}^{n_k}|c_j|N_j^{-1}\Bigr)
\le C_q\frac{N_{n_{k-1}}^{1/2}}{\sqrt{n_k-n_{k-1}}}\norm{(c_j)}{\ell^q}.
\end{aligned}
\end{equation}
Taking $c_j=1$ proves \eqref{eq:carrier-bound}.

For $q\ge2+(q_*-2)/k$, we have $1/q-1/2\le-(2+4k/(q_*-2))^{-1}$, so \eqref{eq:layer-endpoints} gives
\begin{equation}\label{eq:frequency-ratio}
N_{n_{k-1}}^{1/2}(n_k-n_{k-1})^{1/q-1/2}
\le\frac{N_{n_{k-1}}^{1/2}}{2^kA N_{n_{k-1}}^{1/2}}
=2^{-k}A^{-1}.
\end{equation}
Since $2+(q_*-2)/k\downarrow2$, for each $q>2$,
\begin{equation}\label{eq:all-q-sums}
\sum_{k\ge1}
N_{n_{k-1}}^{1/2}(n_k-n_{k-1})^{1/q-1/2}<\infty.
\end{equation}
Combining this summability with \eqref{eq:carrier-bound} yields
\[
\sum_{k\ge0}\norm{v_k^0}{\dot B^{-1}_{\infty,q}}
\le C_qA^{-1}+C_q\sum_{k\ge1}
N_{n_{k-1}}^{1/2}(n_k-n_{k-1})^{1/q-1/2}<\infty,
\]
which proves convergence of \eqref{eq:V0}. Taking $q=q_*$ in \eqref{eq:frequency-ratio}, we obtain
\[
\norm{V^0}{\dot B^{-1}_{\infty,q_*}}
\le\sum_{k\ge0}\norm{v_k^0}{\dot B^{-1}_{\infty,q_*}}
\le C_{q_*}A^{-1}\left(1+\sum_{k\ge1}2^{-k}\right)
\le C_{q_*}A^{-1}.\qedhere
\]
\end{proof}

\subsection{Principal parts}

For $t\ge0$, set $v_0(t)=e^{-(2\pi N_0)^2t}v_0^0$. For $k\ge1$ and $n_{k-1}<j\le n_k$, define
\begin{align}
v_{k,j}(t)&=e^{-(2\pi N_j)^2t}v_{k,j}^0,
\qquad
v_k(t)=\sum_{j=n_{k-1}+1}^{n_k}v_{k,j}(t),
\label{eq:flows}\\
\bar v_{k-1}(t)&=\frac1{n_k-n_{k-1}}
\sum_{j=n_{k-1}+1}^{n_k}e^{-2(2\pi N_j)^2t}v_{k-1}^0.
\label{eq:averaged-flow}
\end{align}
Consider the series
\begin{equation}\label{eq:principal}
v^{(1)}=\sum_{\substack{k\ge0\\k\ \mathrm{even}}}v_k
+\sum_{\substack{k\ge1\\k\ \mathrm{odd}}}\bar v_k,
\qquad
v^{(2)}=\sum_{\substack{k\ge1\\k\ \mathrm{odd}}}v_k
+\sum_{\substack{k\ge0\\k\ \mathrm{even}}}\bar v_k.
\end{equation}

\begin{proposition}[Principal parts]\label{prop:construction}
The series \eqref{eq:principal} converge to real-valued, mean-zero, divergence-free fields. For $i=1,2$ and $2<q\le\infty$,
\[
v^{(i)}\in C^\infty((0,1]\times\T^2)
\cap C([0,1];\dot B^{-1}_{\infty,q})
\cap L^{2,q}(0,1;L^\infty),
\]
with $v^{(i)}(0)=V^0$. Moreover,
\begin{equation}\label{eq:construction-main}
\norm{v^{(i)}}{L^{2,q_*}(0,1;L^\infty)}\le C_{q_*}A^{-1},
\qquad
\sup_{0<t\le1}\sqrt t\,\norm{v^{(i)}(t)}{L^\infty}\le CA^{-1}.
\end{equation}
For $t_0=(2\pi A)^{-2}$ and sufficiently large $A$,
\begin{equation}\label{eq:separation-main}
\norm{v^{(1)}(t_0)-v^{(2)}(t_0)}{L^\infty}\ge e^{-1}.
\end{equation}
\end{proposition}

\begin{proof}
The fields are real-valued, mean-zero and divergence-free by Proposition~\ref{prop:initial}. For $k\ge1$ and $2<q\le\infty$, \eqref{eq:carrier-multipliers} and \eqref{eq:averaged-flow} give
\[
\sup_{0\le t\le1}\left(
\norm{v_k(t)}{\dot B^{-1}_{\infty,q}}
+\norm{\bar v_k(t)}{\dot B^{-1}_{\infty,q}}\right)
\le C_qN_{n_{k-1}}^{1/2}(n_k-n_{k-1})^{1/q-1/2}.
\]
Summing in $k$ and using \eqref{eq:all-q-sums} gives
\[
\sum_{k\ge0}\left(
\norm{v_k}{C([0,1];\dot B^{-1}_{\infty,q})}
+\norm{\bar v_k}{C([0,1];\dot B^{-1}_{\infty,q})}\right)<\infty.
\]
The two series in \eqref{eq:principal} therefore converge in $C([0,1];\dot B^{-1}_{\infty,q})$. Since $v_k(0)=\bar v_k(0)=v_k^0$,
\[
v^{(1)}(0)=v^{(2)}(0)=\sum_{k\ge0}v_k^0=V^0.
\]

For integers $s\ge2$, the product and composition estimates, together with \eqref{eq:amplitude-induction} and \eqref{eq:inverse-bounds}, give
\[
\norm{v_{k,j}^0}{H^s}
\le C_s\frac{N_{n_{k-1}}^{1/2}}{\sqrt{n_k-n_{k-1}}}
\left(N_j^{s+1}+\frac{\norm{v_{k-1}^0}{H^s}}{N_{n_{k-1}}}\right).
\]
Since
\[
\norm{v_0^0}{H^s}\le C_sN_0^s,
\qquad
\frac{N_{n_{k-1}}^{1/2}}{\sqrt{n_k-n_{k-1}}}\le1,
\qquad
\sum_{j=n_{k-1}+1}^{n_k}N_j^{s+1}\le2N_{n_k}^{s+1},
\]
induction gives, for sufficiently large $C_s\ge1$,
\[
\norm{v_{k,j}^0}{H^s}\le C_s^{k+1}N_j^{s+1},
\qquad \norm{v_k^0}{H^s}\le C_s^{k+1}N_{n_k}^{s+1}.
\]
Differentiating \eqref{eq:flows}--\eqref{eq:averaged-flow}, we obtain, for integers $r\ge0$ and $0<\delta<1$,
\[
\begin{aligned}
&\sum_{k\ge1}\sup_{\delta\le t\le1}
\left(\norm{\partial_t^rv_k(t)}{H^s}
+\norm{\partial_t^r\bar v_{k-1}(t)}{H^s}\right)\\
&\quad\le C_{r,s}\sum_{k\ge1}C_s^{k+1}
\sum_{j=n_{k-1}+1}^{n_k}N_j^{s+2r+1}e^{-cN_j^2\delta}
\le C_{r,s}\sum_{j\ge1}C_s^{j+1}N_j^{s+2r+1}e^{-cN_j^2\delta}<\infty.
\end{aligned}
\]
Thus $v^{(i)}\in C^\infty((0,1]\times\T^2)$.

For $k\ge1$, \eqref{eq:carrier-infinity}, \eqref{eq:frequency-ratio} and \eqref{eq:averaged-flow} give
\begin{align*}
\norm{v_k(t)}{L^\infty}
&\le C\frac{N_{n_{k-1}}^{1/2}}{\sqrt{n_k-n_{k-1}}}
\sum_{j=n_{k-1}+1}^{n_k}N_je^{-(2\pi N_j)^2t},\\[3pt]
\norm{\bar v_k(t)}{L^\infty}
&\le\norm{v_k^0}{L^\infty}e^{-2(2\pi N_{n_k+1})^2t}
\\*
&\le C2^{-k}A^{-1}N_{n_k}e^{-2(2\pi N_{n_k+1})^2t}.
\end{align*}
By \eqref{eq:lacunary-heat},
\begin{align*}
\norm{v_k}{L^{2,q}(0,1;L^\infty)}
&\le C_qN_{n_{k-1}}^{1/2}(n_k-n_{k-1})^{1/q-1/2},\\[3pt]
\norm{\bar v_k}{L^{2,q}(0,1;L^\infty)}
&\le C_q2^{-k}A^{-1}\frac{N_{n_k}}{N_{n_k+1}}.
\end{align*}
For the seed flows,
\[
\norm{v_0}{L^{2,q}(0,1;L^\infty)}
+\norm{\bar v_0}{L^{2,q}(0,1;L^\infty)}
\le C_q(N_0^{-1}+N_1^{-1})\le C_qA^{-1}.
\]
Combining these bounds with \eqref{eq:all-q-sums}, we obtain
\[
\norm{v^{(i)}}{L^{2,q}(0,1;L^\infty)}
\le C_q\left[A^{-1}+\sum_{k\ge1}
N_{n_{k-1}}^{1/2}(n_k-n_{k-1})^{1/q-1/2}\right]<\infty.
\]
By \eqref{eq:frequency-ratio} and \eqref{eq:lacunary-heat},
\begin{align*}
\norm{v^{(i)}}{L^{2,q_*}(0,1;L^\infty)}
&\le C_{q_*}A^{-1}\left(1+\sum_{k\ge1}2^{-k}\right)\le C_{q_*}A^{-1},\\[3pt]
\sup_{0<t\le1}\sqrt t\,\norm{v^{(i)}(t)}{L^\infty}
&\le CA^{-1}+C\sum_{k\ge1}\frac{N_{n_{k-1}}^{1/2}}{\sqrt{n_k-n_{k-1}}}
\le CA^{-1}.
\end{align*}

Finally, at $t_0=(2\pi A)^{-2}$, Parseval's identity and \eqref{eq:geometry} give
\[
\norm{v_0(t_0)}{L^\infty}\ge\norm{v_0(t_0)}{L^2}
=e^{-1}\left(\sum_{\xi\in\Lambda}a_\xi(\Id)^2\right)^{1/2}
=\sqrt2\,e^{-1}.
\]
By \eqref{eq:lacunary-heat}, the remaining flows satisfy
\[
\begin{aligned}
&\sum_{k\ge1}\norm{v_k(t_0)}{L^\infty}
+\sum_{k\ge0}\norm{\bar v_k(t_0)}{L^\infty}\\
&\quad\le CA^{-1}\sum_{j\ge1}N_je^{-(2\pi N_j)^2t_0}
+Ce^{-2(2\pi N_1)^2t_0}\\
&\quad\le Ce^{-\frac12(N_1/N_0)^2}
\left(A^{-1}\sum_{j\ge1}N_je^{-(2\pi N_j)^2t_0/2}+1\right)
\le Ce^{-\frac12(N_1/N_0)^2}.
\end{aligned}
\]
Hence, for sufficiently large $A$,
\begin{align*}
&\norm{v^{(1)}(t_0)-v^{(2)}(t_0)}{L^\infty}\\*
&\quad\ge\norm{v_0(t_0)}{L^\infty}
-\sum_{k\ge1}\norm{v_k(t_0)}{L^\infty}
-\sum_{k\ge0}\norm{\bar v_k(t_0)}{L^\infty}\\*
&\quad\ge\sqrt2\,e^{-1}-Ce^{-\frac12(N_1/N_0)^2}\ge e^{-1}.\qedhere
\end{align*}
\end{proof}

\subsection{Residual estimates}

For $0<T\le1$, let $X_T$ be the Banach space of real-valued fields in $C((0,T];L^\infty)$ with norm
\begin{equation*}
\norm\omega{X_T}=\sup_{0<t\le T}t^{1/4}\norm{\omega(t)}{L^\infty}.
\end{equation*}
The residuals of \eqref{eq:principal} are
\begin{equation*}
\mathcal F^{(i)}=-\left[\partial_tv^{(i)}-\Delta v^{(i)}
+\Pj\diver(v^{(i)}\otimes v^{(i)})\right].
\end{equation*}
\begin{proposition}[Residual estimates]\label{prop:residual}
For $i=1,2$ and sufficiently large $A$,
\begin{equation}\label{eq:residual-small}
\norm{\int_0^t e^{(t-s)\Delta}\mathcal F^{(i)}(s)\,ds}{X_1}
+\norm{\mathcal F^{(i)}}{L^1(0,1;\dot B^{-1}_{\infty,2})}\le A^{-1}.
\end{equation}
Moreover, for $2<q\le\infty$,
\begin{equation}\label{eq:principal-nonlinear-integrable}
\Pj\diver(v^{(i)}\otimes v^{(i)})
\in L^1(0,1;\dot B^{-1}_{\infty,q}).
\end{equation}
\end{proposition}

\begin{proof}
All series below converge smoothly for $t>0$ by Proposition~\ref{prop:construction}. For $i=1$, set
\begin{align*}
\mathcal F_1^{(1)}
&=\sum_{\substack{k\ge1\\k\ \mathrm{odd}}}\Delta\bar v_k,
\\
\mathcal F_2^{(1)}
&=-\sum_{\substack{k\ge2\\k\ \mathrm{even}}}
\left[\partial_t\bar v_{k-1}
+\Pj\diver\left(\sum_{j=n_{k-1}+1}^{n_k}v_{k,j}\otimes v_{k,j}\right)\right],
\\
\mathcal F_3^{(1)}
&=-\Pj\diver\Biggl(v^{(1)}\otimes v^{(1)}-v_0\otimes v_0
-\sum_{\substack{k\ge2\\k\ \mathrm{even}}}
\sum_{j=n_{k-1}+1}^{n_k}v_{k,j}\otimes v_{k,j}\Biggr),
\\
\mathcal F_4^{(1)}
&=-\sum_{\substack{k\ge2\\k\ \mathrm{even}}}\sum_{j=n_{k-1}+1}^{n_k}(\partial_t-\Delta)v_{k,j}.
\end{align*}
For $i=2$, replace even $k\ge2$ by odd $k\ge1$ in $\mathcal F_2^{(1)}$ and $\mathcal F_4^{(1)}$, and set
\[
\begin{aligned}
\mathcal F_1^{(2)}&=\sum_{\substack{k\ge0\\k\ \mathrm{even}}}\Delta\bar v_k,\\
\mathcal F_3^{(2)}&=-\Pj\diver\left(v^{(2)}\otimes v^{(2)}
-\sum_{\substack{k\ge1\\k\ \mathrm{odd}}}
\sum_{j=n_{k-1}+1}^{n_k}v_{k,j}\otimes v_{k,j}\right).
\end{aligned}
\]
By \eqref{eq:seed} and Lemma~\ref{lem:planar-cancellation},
\[
(\partial_t-\Delta)v_0=0,\qquad \Pj\diver(v_0\otimes v_0)=0,
\qquad \mathcal F^{(i)}=\sum_{\ell=1}^4\mathcal F_\ell^{(i)}.
\]
We prove the estimates for $i=1$; the case $i=2$ is analogous.

For $\mathcal F_1^{(1)}$, \eqref{eq:averaged-flow} gives
\[
\mathcal F_1^{(1)}(t)
=\sum_{\substack{k\ge1\\k\ \mathrm{odd}}}
\frac1{n_{k+1}-n_k}\sum_{j=n_k+1}^{n_{k+1}}
e^{-2(2\pi N_j)^2t}\Delta v_k^0.
\]
For $k\ge0$, \eqref{eq:seed}, \eqref{eq:spatial-bounds} and \eqref{eq:fast-heat} give
\begin{align*}
&\norm{\int_0^t e^{(t-s)\Delta}\Delta\bar v_k(s)\,ds}{X_1}
+\norm{\Delta\bar v_k}{L^1(0,1;\dot B^{-1}_{\infty,2})}\\*
&\quad\le C\int_0^1\norm{\Delta\bar v_k(t)}{L^\infty}\,dt
\\*
&\quad\le C N_{n_k}^{m+1}\int_0^1e^{-2(2\pi N_{n_k+1})^2t}\,dt
\le C N_{n_k}^{m+1}N_{n_k+1}^{-2}.
\end{align*}
Summing over odd $k\ge1$, we obtain
\[
\begin{aligned}
&\norm{\int_0^t e^{(t-s)\Delta}\mathcal F_1^{(1)}(s)\,ds}{X_1}
+\norm{\mathcal F_1^{(1)}}{L^1(0,1;\dot B^{-1}_{\infty,2})}\\
&\quad\le C\sum_{k\ge0}N_{n_k}^{m+1}N_{n_k+1}^{-2}.
\end{aligned}
\]

For $\mathcal F_2^{(1)}$, \eqref{eq:averaged-flow} gives
\[
\mathcal F_2^{(1)}=\sum_{\substack{k\ge2\\k\ \mathrm{even}}}
\sum_{j=n_{k-1}+1}^{n_k}\mathcal F_{k,j},
\]
where, for even $k\ge2$ and $n_{k-1}<j\le n_k$,
\begin{equation}\label{eq:cancel}
\begin{aligned}
\mathcal F_{k,j}(t)
&=e^{-2(2\pi N_j)^2t}\left[
\frac{2(2\pi N_j)^2}{n_k-n_{k-1}}v_{k-1}^0
-\Pj\diver(v_{k,j}^0\otimes v_{k,j}^0)\right]\\
&=-e^{-2(2\pi N_j)^2t}\Pj\diver\left[
v_{k,j}^0\otimes v_{k,j}^0
-\frac{2(2\pi N_j)^2}{n_k-n_{k-1}}
\left(N_{n_{k-1}}\Id+S(\Delta^{-1}v_{k-1}^0)\right)\right].
\end{aligned}
\end{equation}

Fix $0<c<1/4$ with $4c\le\min_{\xi+\xi'\ne0}|\xi+\xi'|$.
For $R=\Id+S(\Delta^{-1}v_{k-1}^0)/N_{n_{k-1}}$, \eqref{eq:amplitude-induction} and \eqref{eq:fourier-tail} give, for $\xi\in\Lambda$ and $0\le r\le3$,
\begin{equation}\label{eq:amplitude-moments}
\begin{aligned}
\sum_h(1+|h|)^r|\widehat{a_\xi(R)}(h)|
&\le C N_{j-1}^m,\\
\sum_{|h|\ge cN_j}(1+|h|)^r|\widehat{a_\xi(R)}(h)|
&\le C N_{j-1}^mN_j^{r+1-m}.
\end{aligned}
\end{equation}
The identities \eqref{eq:opposite-pair}--\eqref{eq:projected-pair} yield
\[
\begin{aligned}
&\sum_{\xi\in\Lambda}
\nabla^\perp\!\left(a_\xi(R)e^{2\pi iN_j\xi\cdot x}\right)
\otimes\nabla^\perp\!\left(a_\xi(R)e^{-2\pi iN_j\xi\cdot x}\right)\\
&\quad=(2\pi N_j)^2R
+\sum_{\xi\in\Lambda}\bigl[\nabla^\perp a_\xi(R)\bigr]^{\otimes2},
\end{aligned}
\]
and, for $\xi+\xi'\ne0$,
\[
\begin{aligned}
&\Pj\diver\left[
a_\xi(R)a_{\xi'}(R)
(\bar\xi\otimes\bar\xi'+\bar\xi'\otimes\bar\xi)
e^{2\pi iN_j(\xi+\xi')\cdot x}\right]\\
&\quad=\Pj\Bigl[\bigl(\bar\xi\otimes\bar\xi'+\bar\xi'\otimes\bar\xi
-(\xi\cdot\xi'-1)\Id\bigr)
\nabla\bigl(a_\xi(R)a_{\xi'}(R)\bigr)e^{2\pi iN_j(\xi+\xi')\cdot x}\Bigr].
\end{aligned}
\]
Substituting these identities into \eqref{eq:cancel} and using $N_{n_{k-1}}/(n_k-n_{k-1})\le1$, we obtain
\[
\begin{aligned}
&e^{2(2\pi N_j)^2t}|\widehat{\mathcal F_{k,j}}(t,n)|\\
&\quad\le C\sum_{\xi\in\Lambda}\sum_{h+h'=n}
(|h|+|h'|)^3|\widehat{a_\xi(R)}(h)\widehat{a_\xi(R)}(h')|\\
&\qquad+C\sum_{\substack{\xi,\xi'\in\Lambda\\\xi+\xi'\ne0}}
\sum_{h+h'=n-N_j(\xi+\xi')}\sum_{r=1}^3
N_j^{3-r}(|h|+|h'|)^r
|\widehat{a_\xi(R)}(h)\widehat{a_{\xi'}(R)}(h')|.
\end{aligned}
\]
For $\xi+\xi'\ne0$, $|N_j(\xi+\xi')+h+h'|\notin[cN_j,4N_j]$
implies $\max\{|h|,|h'|\}>cN_j$.
By \eqref{eq:amplitude-moments}, for $r=1,2,3$,
\[
\begin{aligned}
\sum_{h,h'}(|h|+|h'|)^r
|\widehat{a_\xi(R)}(h)\widehat{a_{\xi'}(R)}(h')|
&\le C N_{j-1}^{2m},\\
\sum_{\max\{|h|,|h'|\}>cN_j}(|h|+|h'|)^r
|\widehat{a_\xi(R)}(h)\widehat{a_{\xi'}(R)}(h')|
&\le C N_{j-1}^{2m}N_j^{r+1-m}.
\end{aligned}
\]
Combining the Fourier bound with these moment estimates gives
\begin{align*}
\sum_{cN_j\le|n|\le4N_j}|\widehat{\mathcal F_{k,j}}(t,n)|
&\le C N_{j-1}^{2m}\left(1+\sum_{r=1}^3N_j^{3-r}\right)e^{-cN_j^2t}\\*
&\le C N_{j-1}^{2m}N_j^2e^{-cN_j^2t},\\[3pt]
\sum_{|n|\notin[cN_j,4N_j]}|\widehat{\mathcal F_{k,j}}(t,n)|
&\le C N_{j-1}^{2m}\left(1+\sum_{r=1}^3N_j^{3-r}N_j^{r+1-m}\right)e^{-cN_j^2t}\\*
&\le C N_{j-1}^{2m}(1+N_j^{4-m})e^{-cN_j^2t}.
\end{align*}
Applying \eqref{eq:high-heat} to the annular part and \eqref{eq:fast-heat} to its complement gives
\begin{align*}
&\norm{\int_0^t e^{(t-s)\Delta}\mathcal F_{k,j}(s)\,ds}{X_1}
+\norm{\mathcal F_{k,j}}{L^1(0,1;\dot B^{-1}_{\infty,2})}\\*
&\quad\le C N_{j-1}^{2m}
\left[N_j^2\bigl(N_j^{-5/2}+N_j^{-3}\bigr)
+\bigl(1+N_j^{4-m}\bigr)N_j^{-2}\right]
\le C N_{j-1}^{2m}N_j^{-1/2}.
\end{align*}
Summing over $j$ and even $k\ge2$, we obtain
\begin{align*}
&\norm{\int_0^t e^{(t-s)\Delta}\mathcal F_2^{(1)}(s)\,ds}{X_1}
+\norm{\mathcal F_2^{(1)}}{L^1(0,1;\dot B^{-1}_{\infty,2})}\\*
&\quad\le C\sum_{\substack{k\ge2\\k\ \mathrm{even}}}
\sum_{j=n_{k-1}+1}^{n_k}N_{j-1}^{2m}N_j^{-1/2}
\le C\sum_{j\ge1}N_{j-1}^{2m}N_j^{-1/2}.
\end{align*}

For $\mathcal F_3^{(1)}$, expanding \eqref{eq:principal} gives
\[
\begin{aligned}
\mathcal F_3^{(1)}
={}&-\Pj\diver\Biggl[
\sum_{\substack{k\ge2\\k\ \mathrm{even}}}
\sum_{n_{k-1}<j'<j\le n_k}
\bigl(v_{k,j}\otimes v_{k,j'}+v_{k,j'}\otimes v_{k,j}\bigr)\\
&\quad+\sum_{\substack{0\le k'<k\\k,k'\ \mathrm{even}}}
\bigl(v_k\otimes v_{k'}+v_{k'}\otimes v_k\bigr)\\
&\quad+\sum_{\substack{k\ge1\ \mathrm{odd}\\k'\ge0\ \mathrm{even}}}
\bigl(\bar v_k\otimes v_{k'}+v_{k'}\otimes\bar v_k\bigr)
+\sum_{\substack{k,k'\ge1\\k,k'\ \mathrm{odd}}}
\bar v_k\otimes\bar v_{k'}\Biggr].
\end{aligned}
\]
In the first two sums, expand $v_k=\sum_{j=n_{k-1}+1}^{n_k}v_{k,j}$ for $k\ge1$. For $k,k'\ge1$, $n_{k-1}<j\le n_k$, $n_{k'-1}<j'\le n_{k'}$ and $j'<j$, \eqref{eq:amplitude-induction}--\eqref{eq:spatial-bounds} give
\[
\norm{\nabla^r(v_{k,j}\otimes v_{k',j'})(t)}{L^\infty}
\le C N_{j'}^{m+1}N_j^{r+1}e^{-(2\pi N_j)^2t},
\qquad r=0,1.
\]
Splitting the time integral at $t/2$ and applying \eqref{eq:heat-div}, we obtain
\begin{align*}
&\norm{\int_0^t e^{(t-s)\Delta}\Pj\diver(v_{k,j}\otimes v_{k',j'})(s)\,ds}{X_1}\\*
&\quad\le C N_{j'}^{m+1}N_j
\sup_{0<t\le1}t^{1/4}\int_0^t(t-s)^{-1/2}e^{-cN_j^2s}\,ds\\*
&\quad\le C N_{j'}^{m+1}N_j\sup_{0<t\le1}
\left[t^{-1/4}\int_0^{t/2}e^{-cN_j^2s}\,ds
+t^{3/4}e^{-cN_j^2t/2}\right]\\*
&\quad\le C N_{j'}^{m+1}N_j\sup_{0<t\le1}
\min\{t^{3/4},N_j^{-2}t^{-1/4}\}
\le C N_{j'}^{m+1}N_j^{-1/2}.
\end{align*}
Applying \eqref{eq:div-interpolation} to the same product yields
\[
\begin{aligned}
&\norm{\Pj\diver(v_{k,j}\otimes v_{k',j'})}{L^1(0,1;\dot B^{-1}_{\infty,2})}\\
&\quad\le C\int_0^1
\norm{(v_{k,j}\otimes v_{k',j'})(t)}{L^\infty}^{1/2}
\norm{(v_{k,j}\otimes v_{k',j'})(t)}{W^{1,\infty}}^{1/2}\,dt\\
&\quad\le C N_{j'}^{m+1}N_j^{3/2}\int_0^1e^{-cN_j^2t}\,dt
\le C N_{j'}^{m+1}N_j^{-1/2}.
\end{aligned}
\]
The same bounds hold for $v_{k,j}\otimes v_0$ and $v_{k,j}\otimes\bar v_{k'}$, $k'<k$, with $N_{j'}$ replaced by $N_0$ and $N_{n_{k'}}$, respectively. They also hold with the factors interchanged.

For $\bar v_k\otimes v_{k'}$ and $\bar v_k\otimes\bar v_{k'}$ with $k'\le k$, applying \eqref{eq:spatial-bounds} to \eqref{eq:flows}--\eqref{eq:averaged-flow} yields
\begin{equation}\label{eq:cascade-error}
\begin{aligned}
&\sum_{k'=0}^k\Bigl(
\norm{\Pj\diver(\bar v_k\otimes v_{k'})(t)}{L^\infty}
+\norm{\Pj\diver(\bar v_k\otimes\bar v_{k'})(t)}{L^\infty}\Bigr)\\
&\quad\le C\norm{\bar v_k(t)}{H^m}
\sum_{k'=0}^k\bigl(\norm{v_{k'}(t)}{H^m}+\norm{\bar v_{k'}(t)}{H^m}\bigr)\\
&\quad\le C N_{n_k}^{2m+2}e^{-2(2\pi N_{n_k+1})^2t}.
\end{aligned}
\end{equation}
Applying \eqref{eq:fast-heat} to the products in \eqref{eq:cascade-error} and summing all the preceding bounds gives
\begin{align*}
&\norm{\int_0^t e^{(t-s)\Delta}\mathcal F_3^{(1)}(s)\,ds}{X_1}
+\norm{\mathcal F_3^{(1)}}{L^1(0,1;\dot B^{-1}_{\infty,2})}\\*
&\quad\le C\sum_{k\ge0}N_{n_k}^{2m+2}
\int_0^1e^{-2(2\pi N_{n_k+1})^2t}\,dt\\*
&\qquad+C\sum_{j\ge1}N_j^{-1/2}
\Bigl(N_0^{m+1}+\sum_{j'=1}^{j-1}N_{j'}^{m+1}
+\sum_{\substack{k'\ge1\\n_{k'}<j}}N_{n_{k'}}^{m+1}\Bigr)\\*
&\quad\le C\sum_{k\ge0}N_{n_k}^{2m+2}N_{n_k+1}^{-2}
+C\sum_{j\ge1}N_{j-1}^{m+1}N_j^{-1/2}.
\end{align*}

For $\mathcal F_4^{(1)}$, \eqref{eq:flows} gives
\[
\begin{aligned}
\mathcal F_4^{(1)}(t)
&=\sum_{\substack{k\ge2\\k\ \mathrm{even}}}
\sum_{j=n_{k-1}+1}^{n_k}e^{-(2\pi N_j)^2t}
\bigl[\Delta+(2\pi N_j)^2\bigr]v_{k,j}^0,\\
\widehat{(\partial_t-\Delta)v_{k,j}}(t,n)
&=4\pi^2\bigl(|n|^2-N_j^2\bigr)e^{-(2\pi N_j)^2t}\widehat{v_{k,j}^0}(n).
\end{aligned}
\]
Since $|\xi|=1$,
\[
\begin{aligned}
\bigl|N_j^2-|N_j\xi+h|^2\bigr|\,|N_j\xi+h|
&=\bigl|2N_j\xi\cdot h+|h|^2\bigr|\,|N_j\xi+h|\\
&\le C\bigl(N_j^2|h|+N_j|h|^2+|h|^3\bigr).
\end{aligned}
\]
Since $|n-N_j\xi|>cN_j$ for $|n|\notin[cN_j,4N_j]$,
\eqref{eq:carrier} and \eqref{eq:amplitude-moments} give
\begin{align*}
\sum_{cN_j\le|n|\le4N_j}\bigl|N_j^2-|n|^2\bigr|\,|\widehat{v_{k,j}^0}(n)|
&\le C N_{j-1}^m(N_j^2+N_j+1)
\le C N_{j-1}^mN_j^2,\\[3pt]
\sum_{|n|\notin[cN_j,4N_j]}\bigl|N_j^2-|n|^2\bigr|\,|\widehat{v_{k,j}^0}(n)|
&\le C N_{j-1}^m\sum_{r=1}^3N_j^{3-r}N_j^{r+1-m}
\le C N_{j-1}^mN_j^{4-m}.
\end{align*}

By \eqref{eq:high-heat} and \eqref{eq:fast-heat},
\begin{align*}
&\norm{\int_0^t e^{(t-s)\Delta}(\partial_s-\Delta)v_{k,j}(s)\,ds}{X_1}
+\norm{(\partial_t-\Delta)v_{k,j}}{L^1(0,1;\dot B^{-1}_{\infty,2})}\\*
&\quad\le C N_{j-1}^m
\left[N_j^2\bigl(N_j^{-5/2}+N_j^{-3}\bigr)
+N_j^{4-m}N_j^{-2}\right]
\le C N_{j-1}^mN_j^{-1/2}.
\end{align*}
Summing over $j$ and even $k\ge2$ gives
\begin{align*}
&\norm{\int_0^t e^{(t-s)\Delta}\mathcal F_4^{(1)}(s)\,ds}{X_1}
+\norm{\mathcal F_4^{(1)}}{L^1(0,1;\dot B^{-1}_{\infty,2})}\\*
&\quad\le C\sum_{\substack{k\ge2\\k\ \mathrm{even}}}
\sum_{j=n_{k-1}+1}^{n_k}N_{j-1}^mN_j^{-1/2}
\le C\sum_{j\ge1}N_{j-1}^mN_j^{-1/2}.
\end{align*}

Summing the four bounds, with $m=14$ and $N_j=N_{j-1}^{320}$, gives
\begin{align*}
&\norm{\int_0^t e^{(t-s)\Delta}\mathcal F^{(1)}(s)\,ds}{X_1}
+\norm{\mathcal F^{(1)}}{L^1(0,1;\dot B^{-1}_{\infty,2})}\\*
&\quad\le C\sum_{j\ge1}N_{j-1}^{2m+2}N_j^{-1/2}
\le C\sum_{j\ge1}N_j^{-1/4}
\le CA^{-80}\le A^{-1}.
\end{align*}

To prove \eqref{eq:principal-nonlinear-integrable}, recall that
\[
\begin{aligned}
\Pj\diver(v^{(1)}\otimes v^{(1)})
&=-\mathcal F_2^{(1)}-\mathcal F_3^{(1)}
-\sum_{\substack{k\ge1\\k\ \mathrm{odd}}}\partial_t\bar v_k,\\
\Pj\diver(v^{(2)}\otimes v^{(2)})
&=-\mathcal F_2^{(2)}-\mathcal F_3^{(2)}
-\sum_{\substack{k\ge0\\k\ \mathrm{even}}}\partial_t\bar v_k.
\end{aligned}
\]
For $2<q\le\infty$, \eqref{eq:averaged-flow} and \eqref{eq:all-q-sums} give
\[
\begin{aligned}
&\sum_{k\ge0}\int_0^1\norm{\partial_t\bar v_k(t)}{\dot B^{-1}_{\infty,q}}\,dt\\
&\quad=\sum_{k\ge0}\frac{\norm{v_k^0}{\dot B^{-1}_{\infty,q}}}{n_{k+1}-n_k}
\sum_{j=n_k+1}^{n_{k+1}}\bigl(1-e^{-2(2\pi N_j)^2}\bigr)\\
&\quad\le\sum_{k\ge0}\norm{v_k^0}{\dot B^{-1}_{\infty,q}}<\infty.
\end{aligned}
\]
Thus, by $\dot B^{-1}_{\infty,2}\hookrightarrow\dot B^{-1}_{\infty,q}$,
\begin{align*}
&\norm{\Pj\diver(v^{(i)}\otimes v^{(i)})}{L^1(0,1;\dot B^{-1}_{\infty,q})}\\*
&\quad\le C_q\sum_{\ell=2}^3\norm{\mathcal F_\ell^{(i)}}{L^1(0,1;\dot B^{-1}_{\infty,2})}
+\sum_{k\ge0}\norm{v_k^0}{\dot B^{-1}_{\infty,q}}<\infty.\qedhere
\end{align*}
\end{proof}

\section{Construction of the perturbations}\label{sec:correction}

For vector fields $u,v$, define
\begin{equation*}
\mathcal B(u,v)(t)
=\int_0^t e^{(t-s)\Delta}\Pj\diver(u\otimes v)(s)\,ds.
\end{equation*}

\begin{proposition}[Perturbation around a smooth solution]\label{prop:correction}
Let $U$ be a smooth, real-valued, mean-zero, divergence-free solution of \eqref{eq:NS} on $[0,T_0]\times\T^d$, with $T_0>0$. There exists $T\in(0,\min\{T_0,1\}]$, depending only on $U$ and $T_0$, such that for every sufficiently large $A$ there are real-valued, mean-zero, divergence-free perturbations $\omega^{(1)}$ and $\omega^{(2)}$ satisfying
\begin{equation}\label{eq:correction-statement}
\norm{\omega^{(i)}}{X_T}\le4A^{-1},\qquad i=1,2,
\end{equation}
and
\[
\omega^{(i)}\in C([0,T];\dot B^{-1}_{\infty,2})
\cap C^\infty((0,T]\times\T^d),\qquad\omega^{(i)}(0)=0.
\]
The fields $u^{(i)}=U+v^{(i)}+\omega^{(i)}$, $i=1,2$, solve \eqref{eq:NS} with initial datum $U(0)+V^0$ and satisfy \eqref{eq:class}, \eqref{eq:mild} and \eqref{eq:lorentz-class}. If $U=0$ on $[0,1]$, one may take $T=1$.
\end{proposition}

\begin{proof}
Fix $i\in\{1,2\}$ and write $\norm U{L^\infty_{t,x}}=\norm U{L^\infty([0,T_0]\times\T^d)}$.
Writing $u^{(i)}=U+v^{(i)}+\omega^{(i)}$ reduces \eqref{eq:NS} to
\[
(\partial_t-\Delta)\omega^{(i)}
+\Pj\diver\bigl(u^{(i)}\otimes u^{(i)}-U\otimes U-v^{(i)}\otimes v^{(i)}\bigr)=\mathcal F^{(i)}.
\]
By Duhamel's formula, we seek a fixed point of
\begin{equation}\label{eq:background-map}
\begin{aligned}
(\mathcal T\omega)(t)
&=\int_0^t e^{(t-s)\Delta}\mathcal F^{(i)}(s)\,ds
-\mathcal B(U,v^{(i)})(t)-\mathcal B(v^{(i)},U)(t)\\
&\quad-\mathcal B(U+v^{(i)},\omega)(t)-\mathcal B(\omega,U+v^{(i)})(t)
-\mathcal B(\omega,\omega)(t).
\end{aligned}
\end{equation}
For $\omega,\widetilde\omega\in X_T$, \eqref{eq:construction-main} and \eqref{eq:heat-div} give
\begin{align*}
\norm{\mathcal B(v^{(i)},\omega)}{X_T}
&\le CA^{-1}\norm\omega{X_T}
\sup_{0<t\le T}t^{1/4}\int_0^t(t-s)^{-1/2}s^{-3/4}\,ds\\*
&\le CA^{-1}\norm\omega{X_T}
\int_0^1(1-s)^{-1/2}s^{-3/4}\,ds
\le CA^{-1}\norm\omega{X_T}.
\end{align*}
Similarly,
\begin{align*}
\norm{\mathcal B(U,v^{(i)})}{X_T}
&\le C\norm U{L^\infty_{t,x}}A^{-1}T^{1/4},\\[3pt]
\norm{\mathcal B(U,\omega)}{X_T}
&\le C\norm U{L^\infty_{t,x}}T^{1/2}\norm\omega{X_T},\\[3pt]
\norm{\mathcal B(\omega,\widetilde\omega)}{X_T}
&\le CT^{1/4}\norm\omega{X_T}\norm{\widetilde\omega}{X_T}.
\end{align*}
The same bounds hold with the factors interchanged, and the heat potentials are continuous for $t>0$. Thus $\mathcal T(X_T)\subset X_T$. Choose $T$ such that
\[
0<T\le\min\{T_0,1\},
\qquad C\norm U{L^\infty_{t,x}}T^{1/4}\le\tfrac18,
\]
and then choose $A$ sufficiently large. For $\norm\omega{X_T},\norm{\widetilde\omega}{X_T}\le4A^{-1}$, \eqref{eq:residual-small} gives
\begin{gather*}
\begin{aligned}
\norm{\mathcal T\omega}{X_T}
&\le A^{-1}+C\norm U{L^\infty_{t,x}}A^{-1}T^{1/4}
+C\bigl(\norm U{L^\infty_{t,x}}T^{1/2}+A^{-1}\bigr)\norm\omega{X_T}\\
&\quad+CT^{1/4}\norm\omega{X_T}^2
\le2A^{-1}+\tfrac12\norm\omega{X_T}\le4A^{-1},
\end{aligned}\\[3pt]
\begin{aligned}
&\norm{\mathcal T\omega-\mathcal T\widetilde\omega}{X_T}\\
&\quad\le C\left[\norm U{L^\infty_{t,x}}T^{1/2}+A^{-1}
+T^{1/4}\bigl(\norm\omega{X_T}+\norm{\widetilde\omega}{X_T}\bigr)\right]
\norm{\omega-\widetilde\omega}{X_T}\\
&\quad\le\tfrac12\norm{\omega-\widetilde\omega}{X_T}.
\end{aligned}
\end{gather*}
Banach's fixed-point theorem gives a mean-zero, divergence-free $\omega^{(i)}=\mathcal T\omega^{(i)}$ with \eqref{eq:correction-statement}. If $U=0$, these bounds hold with $T=1$.

For the initial trace, \eqref{eq:construction-main} and \eqref{eq:correction-statement} give, for $0<s\le T\le1$,
\[
\begin{aligned}
&\norm{\bigl(u^{(i)}\otimes u^{(i)}-U\otimes U-v^{(i)}\otimes v^{(i)}\bigr)(s)}{L^\infty}\\
&\quad\le C\left[\norm U{L^\infty_{t,x}}A^{-1}(s^{-1/2}+s^{-1/4})
+A^{-2}(s^{-3/4}+s^{-1/2})\right]\\
&\quad\le C\left(\norm U{L^\infty_{t,x}}A^{-1}s^{-1/2}+A^{-2}s^{-3/4}\right).
\end{aligned}
\]
For $\alpha=1/2,3/4$, the change of variables $s=tr$ gives
\[
\begin{aligned}
\int_0^t\sqrt{\log(e/(t-s))}\,s^{-\alpha}\,ds
&=t^{1-\alpha}\int_0^1
\sqrt{\log(e/t)+\log(1/(1-r))}\,r^{-\alpha}\,dr\\
&\le Ct^{1-\alpha}\sqrt{\log(e/t)}.
\end{aligned}
\]
Thus \eqref{eq:log-heat} and \eqref{eq:background-map} give
\[
\begin{aligned}
\norm{\omega^{(i)}(t)}{\dot B^{-1}_{\infty,2}}
&\le C\int_0^t\norm{\mathcal F^{(i)}(s)}{\dot B^{-1}_{\infty,2}}\,ds\\
&\quad+C\int_0^t\sqrt{\log(e/(t-s))}
\left(A^{-2}s^{-3/4}+\norm U{L^\infty_{t,x}}A^{-1}s^{-1/2}\right)\,ds\\
&\le C\int_0^t\norm{\mathcal F^{(i)}(s)}{\dot B^{-1}_{\infty,2}}\,ds\\
&\quad+C\sqrt{\log(e/t)}\left(A^{-2}t^{1/4}+\norm U{L^\infty_{t,x}}A^{-1}t^{1/2}\right)
\longrightarrow0\quad(t\downarrow0).
\end{aligned}
\]
Thus $\omega^{(i)}\in C([0,T];\dot B^{-1}_{\infty,2})$ with $\omega^{(i)}(0)=0$.

The semigroup property in \eqref{eq:background-map} gives, for $0<\delta<t\le T$,
\[
u^{(i)}(t)=e^{(t-\delta)\Delta}u^{(i)}(\delta)
-\int_\delta^t e^{(t-s)\Delta}\Pj\diver(u^{(i)}\otimes u^{(i)})(s)\,ds.
\]
Since $u^{(i)}$ is bounded on compact time intervals in $(0,T]$, local smoothing gives $u^{(i)},\omega^{(i)}\in C^\infty((0,T]\times\T^d)$, with
\[
p^{(i)}=-\Delta^{-1}\diver\diver(u^{(i)}\otimes u^{(i)}).
\]
Proposition~\ref{prop:construction} and $\dot B^{-1}_{\infty,2}\hookrightarrow\dot B^{-1}_{\infty,q}$ give \eqref{eq:class}, with
\[
u^{(i)}(0)=U(0)+v^{(i)}(0)+\omega^{(i)}(0)=U(0)+V^0=u_0.
\]

To pass to the initial time in Duhamel's formula, we use \eqref{eq:construction-main}, \eqref{eq:correction-statement} and \eqref{eq:log-heat}:
\begin{align*}
&\int_0^t\norm{e^{(t-s)\Delta}\Pj\diver
\bigl(u^{(i)}\otimes u^{(i)}-v^{(i)}\otimes v^{(i)}\bigr)(s)}{\dot B^{-1}_{\infty,2}}\,ds\\*
&\quad\le C\bigl(1+\norm U{L^\infty_{t,x}}\bigr)^2
\int_0^t\sqrt{\log(e/(t-s))}(1+s^{-3/4})\,ds<\infty.
\end{align*}
Combining this bound with \eqref{eq:principal-nonlinear-integrable} and letting $\delta\downarrow0$ gives \eqref{eq:mild}.

Finally, for $2<q\le\infty$,
\[
\norm{u^{(i)}}{L^{2,q}(0,T;L^\infty)}
\le C_q\left(T^{1/2}\norm U{L^\infty_{t,x}}
+\norm{v^{(i)}}{L^{2,q}(0,T;L^\infty)}+A^{-1}T^{1/4}\right)<\infty.\qedhere
\]
\end{proof}

\section{Proofs of the main results}\label{sec:proofs}

\subsection{Non-uniqueness and density}

\begin{proof}[Proof of Theorem~\ref{thm:background}]
Fix $T\in(0,\min\{T_0,1\}]$ as in Proposition~\ref{prop:correction}, depending only on $U$ and $T_0$.
For $2<q_*<\infty$, that proposition gives two solutions satisfying \eqref{eq:class} and \eqref{eq:mild}, with common initial datum $u_0=U(0)+V^0$.
Given $\varepsilon>0$, choose $A$ so that
\[
C_{q_*}A^{-1}<\varepsilon,
\qquad t_0=(2\pi A)^{-2}<T,
\qquad 8\sqrt{2\pi}\,A^{-1/2}\le\tfrac12e^{-1}.
\]
Then Proposition~\ref{prop:initial} gives
\[
\norm{u_0-U(0)}{\dot B^{-1}_{\infty,q_*}}
=\norm{V^0}{\dot B^{-1}_{\infty,q_*}}\le C_{q_*}A^{-1}<\varepsilon,
\]
and \eqref{eq:separation-main} and \eqref{eq:correction-statement} give
\[
\begin{aligned}
\norm{u^{(1)}(t_0)-u^{(2)}(t_0)}{L^\infty}
&\ge\norm{v^{(1)}(t_0)-v^{(2)}(t_0)}{L^\infty}
-t_0^{-1/4}\sum_{i=1}^2\norm{\omega^{(i)}}{X_T}\\
&\ge e^{-1}-8\sqrt{2\pi}\,A^{-1/2}\ge\tfrac12e^{-1}>0.
\end{aligned}
\]
For $U=0$ on $[0,1]$, Proposition~\ref{prop:correction} allows $T=1$.

For $q_*=\infty$, use exponent $4$ in the construction. With the same $T$ and sufficiently large $A$,
\[
\norm{u_0-U(0)}{\dot B^{-1}_{\infty,\infty}}
\le C\norm{V^0}{\dot B^{-1}_{\infty,4}}\le CA^{-1}<\varepsilon.
\]

For density, let $2<q_*<\infty$ and let $g\in\dot B^{-1}_{\infty,q_*}$ be mean-zero and divergence-free.
Choose a smooth Fourier truncation $f$ with $\norm{f-g}{\dot B^{-1}_{\infty,q_*}}<\varepsilon/2$, and let $U$ be the local classical solution with $U(0)=f$.
Applying the preceding construction to this $U$ with tolerance $\varepsilon/2$ gives $u_0=f+V^0$ and two distinct solutions. Hence
\[
\norm{u_0-g}{\dot B^{-1}_{\infty,q_*}}
\le \norm{f-g}{\dot B^{-1}_{\infty,q_*}}+\norm{V^0}{\dot B^{-1}_{\infty,q_*}}
<\varepsilon.
\]
At $q_*=\infty$ the same argument applies in the closure of the smooth fields. The constructed $u_0$ belongs to this closure since $u_0\in \dot B^{-1}_{\infty,4}$.
\end{proof}

\subsection{Uniqueness}\label{sec:uniqueness}

\begin{proof}[Proof of Theorem~\ref{thm:uniqueness}]
Let $u_0$ be the common initial datum. Positive-time smoothness and the strong initial trace place $u_0$ in the $\dot B^{-1+d/p}_{p,q}$ closure of smooth fields, where spatial translations are continuous. The common mean is conserved, so the Galilean change
\[
u^{(i)}(t,x)\longmapsto u^{(i)}(t,x+ct)-c,
\qquad c=\int_{\T^d}u_0(x)\,dx,
\]
reduces the proof to mean-zero solutions in the same class.

For $0<h\le1/16$, use the local Koch--Tataru norm
\[
\norm z{\mathrm{KT}(h)}
=\sup_{0<t<h}\sqrt t\,\norm{z(t)}{L^\infty}
+\sup_{\substack{x\in\T^d\\0<r\le\sqrt h}}
\left(r^{-d}\int_0^{r^2}\int_{B(x,r)}|z(s,y)|^2\,dy\,ds\right)^{1/2}.
\]
Here $B(x,r)$ is a ball for the flat torus metric. By \cite[Theorem 3 and Lemmas 3.1--3.2]{KT},
\begin{equation}\label{eq:KT-bilinear}
\norm{\mathcal B(w,z)}{\mathrm{KT}(h)}
\le C\norm w{\mathrm{KT}(h)}\norm z{\mathrm{KT}(h)}.
\end{equation}
To obtain the torus estimate, extend the fields periodically to $\mathbb R^d$. The local norms agree for $r\le1/4$, and the periodized differentiated Oseen kernel represents the periodic heat potential.

For mean-zero $f$ and the stated exponents, we first show that
\begin{equation}\label{eq:KT-heat}
\norm{e^{t\Delta}f}{\mathrm{KT}(h)}
\le C\norm f{\dot B^{-1+d/p}_{p,q}}.
\end{equation}
For $p=\infty$ and $q\le2$, dyadic heat estimates and discrete Young's inequality give
\[
\begin{aligned}
\int_0^h\norm{e^{t\Delta}f}{L^\infty}^2\,dt
&\le C\sum_{\ell,\ell'\ge-1}
\frac{\norm{\dot\Delta_\ell f}{L^\infty}\norm{\dot\Delta_{\ell'}f}{L^\infty}}
{2^{2\ell}+2^{2\ell'}}\\
&\le C\sum_{\ell,\ell'\ge-1}2^{-|\ell-\ell'|}
\bigl(2^{-\ell}\norm{\dot\Delta_\ell f}{L^\infty}\bigr)
\bigl(2^{-\ell'}\norm{\dot\Delta_{\ell'}f}{L^\infty}\bigr)\\
&\le C\norm f{\dot B^{-1}_{\infty,2}}^2
\le C\norm f{\dot B^{-1}_{\infty,q}}^2.
\end{aligned}
\]
The heat semigroup is an $L^\infty$ contraction, so
\[
t\norm{e^{t\Delta}f}{L^\infty}^2
\le\int_0^t\norm{e^{s\Delta}f}{L^\infty}^2\,ds,
\qquad
\norm{e^{t\Delta}f}{\mathrm{KT}(h)}
\le C\norm{e^{t\Delta}f}{L^2(0,h;L^\infty)}.
\]

For $p<\infty$, choose $\max\{p,d\}<P<\infty$. Bernstein's inequality gives
\[
\dot B^{-1+d/p}_{p,q}
\hookrightarrow\dot B^{-1+d/P}_{P,\infty}
\hookrightarrow\dot B^{-1}_{\infty,\infty}.
\]
Since $1-d/P>0$, the dyadic heat estimates give
\[
\begin{aligned}
\norm{e^{t\Delta}f}{L^P}
&\le C\sum_{\ell\ge-1}e^{-ct2^{2\ell}}\norm{\dot\Delta_\ell f}{L^P}\\
&\le C\norm f{\dot B^{-1+d/P}_{P,\infty}}
\sum_{\ell\ge-1}2^{(1-d/P)\ell}e^{-ct2^{2\ell}}\\
&\le Ct^{-\frac12+\frac d{2P}}\norm f{\dot B^{-1+d/p}_{p,q}}.
\end{aligned}
\]
Similarly,
\[
\norm{e^{t\Delta}f}{L^\infty}
\le C\norm f{\dot B^{-1}_{\infty,\infty}}
\sum_{\ell\ge-1}2^\ell e^{-ct2^{2\ell}}
\le Ct^{-1/2}\norm f{\dot B^{-1+d/p}_{p,q}}.
\]
Since $P>d\ge2$, H\"older's inequality yields
\[
\begin{aligned}
r^{-d}\int_0^{r^2}\int_{B(x,r)}|e^{s\Delta}f|^2\,dy\,ds
&\le Cr^{-2d/P}\int_0^{r^2}\norm{e^{s\Delta}f}{L^P}^2\,ds\\
&\le Cr^{-2d/P}\norm f{\dot B^{-1+d/p}_{p,q}}^2
\int_0^{r^2}s^{-1+d/P}\,ds
\le C\norm f{\dot B^{-1+d/p}_{p,q}}^2.
\end{aligned}
\]
This proves \eqref{eq:KT-heat}. Applying it to $u_0-u^{(1)}(\delta)$ and the heat contraction to $u^{(1)}(\delta)$ gives, for $\delta>0$,
\[
\norm{e^{t\Delta}u_0}{\mathrm{KT}(h)}
\le C\norm{u_0-u^{(1)}(\delta)}{\dot B^{-1+d/p}_{p,q}}
+C\sqrt h\,\norm{u^{(1)}(\delta)}{L^\infty}.
\]
Letting $h\downarrow0$ gives
\[
\limsup_{h\downarrow0}\norm{e^{t\Delta}u_0}{\mathrm{KT}(h)}
\le C\norm{u_0-u^{(1)}(\delta)}{\dot B^{-1+d/p}_{p,q}}
\longrightarrow0\qquad(\delta\downarrow0).
\]

Let $C\ge1$ dominate the constants in \eqref{eq:KT-bilinear} and \eqref{eq:KT-heat}.
Choose $0<h<\min\{T/2,1/16\}$ so that, for all sufficiently small $0<\delta<h$,
\[
\norm{e^{t\Delta}u^{(i)}(\delta)}{\mathrm{KT}(h)}
\le C\norm{u^{(i)}(\delta)-u_0}{\dot B^{-1+d/p}_{p,q}}
+\norm{e^{t\Delta}u_0}{\mathrm{KT}(h)}\le\frac1{8C},
\qquad i=1,2.
\]
For $\norm w{\mathrm{KT}(h)},\norm z{\mathrm{KT}(h)}\le1/(4C)$, \eqref{eq:KT-bilinear} gives
\begin{align*}
\norm{e^{t\Delta}u^{(i)}(\delta)-\mathcal B(w,w)}{\mathrm{KT}(h)}
&\le\frac1{8C}+C\norm w{\mathrm{KT}(h)}^2
\le\frac3{16C}<\frac1{4C},\\[3pt]
\norm{\mathcal B(w,w)-\mathcal B(z,z)}{\mathrm{KT}(h)}
&\le C\bigl(\norm w{\mathrm{KT}(h)}+\norm z{\mathrm{KT}(h)}\bigr)
\norm{w-z}{\mathrm{KT}(h)}\\*
&\le\tfrac12\norm{w-z}{\mathrm{KT}(h)}.
\end{align*}
By Banach's fixed-point theorem, there is a unique fixed point in this ball with initial datum $u^{(i)}(\delta)$, smooth for positive times. The given solution satisfies
\[
u^{(i)}(\delta+t)=e^{t\Delta}u^{(i)}(\delta)
-\int_0^t e^{(t-s)\Delta}\Pj\diver
\bigl(u^{(i)}(\delta+s)\otimes u^{(i)}(\delta+s)\bigr)\,ds,
\qquad 0<t\le h.
\]
For sufficiently small $s>0$,
\[
\norm{u^{(i)}(\delta+\cdot)}{\mathrm{KT}(s)}
\le C\sqrt s\sup_{\delta\le t\le\delta+s}\norm{u^{(i)}(t)}{L^\infty}
\le\frac1{4C}.
\]
Thus $u^{(i)}(\delta+\cdot)$ coincides with this fixed point for small $s$. Classical uniqueness extends the equality to $(0,h]$, and hence
\[
\norm{u^{(i)}(\delta+\cdot)}{\mathrm{KT}(h)}\le\frac1{4C},
\qquad i=1,2.
\]
Subtracting the Duhamel formulas starting at $\delta$ and applying \eqref{eq:KT-bilinear}, we obtain
\[
\begin{aligned}
\norm{u^{(1)}(\delta+\cdot)-u^{(2)}(\delta+\cdot)}{\mathrm{KT}(h)}
&\le\norm{e^{t\Delta}\bigl(u^{(1)}(\delta)-u^{(2)}(\delta)\bigr)}{\mathrm{KT}(h)}\\
&\quad+\tfrac12\norm{u^{(1)}(\delta+\cdot)-u^{(2)}(\delta+\cdot)}{\mathrm{KT}(h)}.
\end{aligned}
\]
Absorbing the last term and using \eqref{eq:KT-heat}, we obtain, for $0<\delta<t<h$,
\begin{align*}
\sqrt{t-\delta}\,\norm{u^{(1)}(t)-u^{(2)}(t)}{L^\infty}
&\le\norm{u^{(1)}(\delta+\cdot)-u^{(2)}(\delta+\cdot)}{\mathrm{KT}(h)}\\*
&\le2C\norm{u^{(1)}(\delta)-u^{(2)}(\delta)}{\dot B^{-1+d/p}_{p,q}}\\*
&\le2C\sum_{i=1}^2\norm{u^{(i)}(\delta)-u_0}{\dot B^{-1+d/p}_{p,q}}
\longrightarrow0\quad(\delta\downarrow0).
\end{align*}
Thus $u^{(1)}=u^{(2)}$ on $(0,h)$. Classical uniqueness and the common initial trace give equality on $[0,T]$.
\end{proof}

\appendix
\section{Geometric identities}\label{app:geometry}

We use the two-dimensional geometric decomposition and stationary-flow identity recalled in \cite[Lemmas A.1 and A.4]{CLY}; see also \cite[Lemmas 4 and 6]{CDS}.

\begin{lemma}[Geometric decomposition]\label{lem:geometry}
Let
\[
\Lambda=\{\pm(0,-1),\,\pm(1,0),\,\pm(4/5,-3/5),\,\pm(-4/5,-3/5)\},
\qquad \bar\xi=(-\xi_2,\xi_1).
\]
There are smooth positive functions $a_\xi$, with $a_{-\xi}=a_\xi$, satisfying \eqref{eq:geometry} on the symmetric matrices $R$ with $\max_{r,s}|R_{rs}-\Id_{rs}|\le1/10$.
For every smooth symmetric matrix field $R(x)$ in this neighborhood and $N\in5\mathbb N$,
\begin{equation}\label{eq:opposite-pair}
\begin{aligned}
&\sum_{\xi\in\Lambda}
\nabla^\perp\!\left(a_\xi(R)e^{2\pi iN\xi\cdot x}\right)
\otimes\nabla^\perp\!\left(a_\xi(R)e^{-2\pi iN\xi\cdot x}\right)\\
&\quad=(2\pi N)^2R
+\sum_{\xi\in\Lambda}\bigl[\nabla^\perp a_\xi(R)\bigr]^{\otimes2}.
\end{aligned}
\end{equation}
\end{lemma}

\begin{proof}
Take
\[
\begin{alignedat}{2}
a_{\pm(0,-1)}(R)&=\left(\frac{R_{11}}2-\frac9{100}\right)^{1/2},
&\qquad a_{\pm(1,0)}(R)&=\left(\frac{R_{22}}2-\frac4{25}\right)^{1/2},\\
a_{\pm(4/5,-3/5)}(R)&=\left(\frac18+\frac{25}{48}R_{12}\right)^{1/2},
&\qquad a_{\pm(-4/5,-3/5)}(R)&=\left(\frac18-\frac{25}{48}R_{12}\right)^{1/2}.
\end{alignedat}
\]
Here $a_\xi(R)^2\ge7/96>0$. Substitution gives \eqref{eq:geometry}; pairing $\xi$ and $-\xi$ gives \eqref{eq:opposite-pair}.
\end{proof}

\begin{lemma}[Planar cancellation]\label{lem:planar-cancellation}
For $\xi,\xi'\in\Lambda$, $N\in5\mathbb N$ and $f\in C^\infty(\T^2)$,
\begin{equation}\label{eq:projected-pair}
\begin{aligned}
&\Pj\diver\left[f(\bar\xi\otimes\bar\xi'+\bar\xi'\otimes\bar\xi)
e^{2\pi iN(\xi+\xi')\cdot x}\right]\\
&\quad=\Pj\left[\bigl(\bar\xi\otimes\bar\xi'+\bar\xi'\otimes\bar\xi
-(\xi\cdot\xi'-1)\Id\bigr)\nabla f\,e^{2\pi iN(\xi+\xi')\cdot x}\right].
\end{aligned}
\end{equation}
In particular, for real constants $b_\xi=b_{-\xi}$,
\[
\Pj\diver\left[\left(\sum_{\xi\in\Lambda}
i b_\xi\bar\xi e^{2\pi iN\xi\cdot x}\right)^{\otimes2}\right]=0.
\]
\end{lemma}

\begin{proof}
Since $|\xi|=|\xi'|=1$,
\[
(\bar\xi\otimes\bar\xi'+\bar\xi'\otimes\bar\xi)(\xi+\xi')
=(\xi\cdot\xi'-1)(\xi+\xi').
\]
Since $\Pj\nabla=0$,
\[
\Pj\left(f\nabla e^{2\pi iN(\xi+\xi')\cdot x}\right)
=-\Pj\left(e^{2\pi iN(\xi+\xi')\cdot x}\nabla f\right).
\]
The product rule gives \eqref{eq:projected-pair}. Taking $f=b_\xi b_{\xi'}$ and summing proves the last identity.
\end{proof}

\section{Fourier and heat estimates}\label{app:analysis}

For mean-zero vector fields $v$ on $\T^2$ and $s\in\mathbb R$,
\begin{equation}\label{eq:inverse-bounds}
\norm{S(\Delta^{-1}v)}{L^\infty}\le C\norm v{L^\infty},
\qquad
\norm{S(\Delta^{-1}v)}{H^s}\le C_s\norm v{H^{s-1}}.
\end{equation}
The $L^\infty$ bound follows from the integrable kernel of $S\Delta^{-1}$; the $H^s$ bound follows from
\[
|\widehat{S(\Delta^{-1}v)}(n)|\le C|n|^{-1}|\widehat v(n)|,
\qquad n\ne0.
\]

For $f\in H^m(\T^2)$, $N\ge1$, $c>0$ and $0\le r<m-1$, Cauchy--Schwarz gives
\begin{equation}\label{eq:fourier-tail}
\begin{aligned}
\sum_{|h|\ge cN}(1+|h|)^r|\widehat f(h)|
&\le\left(\sum_{|h|\ge cN}(1+|h|)^{-2(m-r)}\right)^{1/2}\norm f{H^m}\\
&\le C_{m,r,c}N^{r+1-m}\norm f{H^m}.
\end{aligned}
\end{equation}
The same argument over all $h\in\mathbb Z^2$ yields
\[
\sum_h(1+|h|)^r|\widehat f(h)|\le C_{m,r}\norm f{H^m}.
\]

For a smooth tensor $G$, dyadic localization gives
\begin{equation}\label{eq:div-interpolation}
\begin{aligned}
\norm{\Pj\diver G}{\dot B^{-1}_{\infty,2}}^2
&\le C\sum_{\ell\ge-1}\min\!\left\{\norm G{L^\infty}^2,
2^{-2\ell}\norm{\nabla G}{L^\infty}^2\right\}\\
&\le C\sum_{\ell\ge-1}2^{-\ell}\norm G{L^\infty}\norm{\nabla G}{L^\infty}
\le C\norm G{L^\infty}\norm G{W^{1,\infty}}.
\end{aligned}
\end{equation}

For a tensor field $G$ on $\T^d$, the differentiated heat kernel estimate gives
\begin{equation}\label{eq:heat-div}
\norm{e^{t\Delta}\Pj\diver G}{L^\infty}
\le Ct^{-1/2}\norm G{L^\infty},\qquad 0<t\le1.
\end{equation}
\begin{lemma}[Heat series]
For $c_j\ge0$ and $1\le q\le\infty$,
\begin{equation}\label{eq:lacunary-heat}
\left\|\sum_{j\ge1}c_jN_je^{-(2\pi N_j)^2t}\right\|_{L^{2,q}(0,1)}
\le C_q\norm{(c_j)}{\ell^q}.
\end{equation}
\end{lemma}

\begin{proof}
Since $N_{j+1}\ge2N_j$,
\[
\sqrt t\sum_{j\ge1}N_je^{-(2\pi N_j)^2t}
\le\sqrt t\sum_{N_j\le t^{-1/2}}N_j
+Ct^{-1/2}\sum_{N_j>t^{-1/2}}N_j^{-1}\le C.
\]
For finitely supported $(c_j)$ and $1\le q<\infty$, monotonicity in $t$ and weighted H\"older's inequality give
\[
\begin{aligned}
\left\|\sum_jc_jN_je^{-(2\pi N_j)^2t}\right\|_{L^{2,q}(0,1)}^q
&\le C_q\int_0^1\left(\sqrt t\sum_jc_jN_je^{-(2\pi N_j)^2t}\right)^q\frac{dt}{t}\\
&\le C_q\sum_jc_j^q\int_0^1N_jt^{-1/2}e^{-(2\pi N_j)^2t}\,dt
\le C_q\sum_jc_j^q.
\end{aligned}
\]
For $q=\infty$,
\[
\sup_{0<t\le1}\sqrt t\sum_jc_jN_je^{-(2\pi N_j)^2t}
\le C\norm{(c_j)}{\ell^\infty}.
\]
The general case follows by monotone convergence.
\end{proof}

\begin{lemma}[Heat potentials]
For every mean-zero $F\in L^1(0,1;L^\infty)$,
\begin{equation}\label{eq:fast-heat}
\norm{\int_0^t e^{(t-s)\Delta}F(s)\,ds}{X_1}
+\norm F{L^1(0,1;\dot B^{-1}_{\infty,2})}
\le C\int_0^1\norm{F(s)}{L^\infty}\,ds.
\end{equation}
If, for $N\ge1$, $B>0$ and fixed $c_0,c_1,C_0>0$,
\[
\supp\widehat F(t,\cdot)\subset\{c_0N\le|n|\le C_0N\},
\qquad \sum_n|\widehat F(t,n)|\le Be^{-c_1N^2t},
\]
then
\begin{equation}\label{eq:high-heat}
\norm{\int_0^t e^{(t-s)\Delta}F(s)\,ds}{X_1}
\le CBN^{-5/2},\qquad
\norm F{L^1(0,1;\dot B^{-1}_{\infty,2})}\le CBN^{-3}.
\end{equation}
\end{lemma}

\begin{proof}
The heat contraction and $L^\infty\hookrightarrow\dot B^{-1}_{\infty,2}$ give
\[
\begin{aligned}
\norm{\int_0^t e^{(t-s)\Delta}F(s)\,ds}{X_1}
&\le\sup_{0<t\le1}t^{1/4}\int_0^t\norm{F(s)}{L^\infty}\,ds
\le\norm F{L^1(0,1;L^\infty)},\\
\norm F{L^1(0,1;\dot B^{-1}_{\infty,2})}
&\le C\norm F{L^1(0,1;L^\infty)}.
\end{aligned}
\]
This proves \eqref{eq:fast-heat}. Under the Fourier support assumption,
\begin{align*}
\norm{\int_0^t e^{(t-s)\Delta}F(s)\,ds}{X_1}
&\le CB\sup_{0<t\le1}t^{1/4}\int_0^t
e^{-cN^2(t-s)}e^{-c_1N^2s}\,ds\\*
&\le CB\sup_{0<t\le1}t^{5/4}e^{-c'N^2t}\le CBN^{-5/2},\\[3pt]
\norm F{L^1(0,1;\dot B^{-1}_{\infty,2})}
&\le CBN^{-1}\int_0^1e^{-c_1N^2t}\,dt\le CBN^{-3}.\qedhere
\end{align*}
\end{proof}

For $0<t\le1$, dyadic localization gives
\begin{equation}\label{eq:log-heat}
\norm{e^{t\Delta}\Pj\diver G}{\dot B^{-1}_{\infty,2}}
\le C\left(\sum_{\ell\ge-1}e^{-ct2^{2\ell}}\right)^{1/2}\norm G{L^\infty}
\le C\sqrt{\log(e/t)}\norm G{L^\infty}.
\end{equation}

\noindent\textbf{Acknowledgments} This work was partially supported by the National Natural Science Foundation of China (No.12571261).

\end{document}